\documentclass[letterpaper, 10 pt, conference]{ieeeconf}  

\IEEEoverridecommandlockouts                              
\usepackage{graphics} 
\usepackage{epsfig} 
\usepackage{times} 
\usepackage{amsmath}

\usepackage{amssymb}
\usepackage{physics}
\usepackage{mathtools}

\usepackage{algorithm}
\usepackage{algpseudocode} 

\usepackage{multicol}

\usepackage{subcaption}
\usepackage{caption}

\usepackage{xcolor}

\newcommand*{\QEDB}{\hfill\ensuremath{\blacksquare}}%

\newcommand*{\vsepfbox}[1]{%
  \begingroup
    \sbox0{\fbox{#1}}%
    \setlength{\fboxrule}{0pt}%
    \mbox{\kern-\fboxsep\fbox{\unhbox0}\kern-\fboxsep}%
  \endgroup
}

\title{\LARGE \bf
BFGS-ADMM for Large-Scale Distributed Optimization
}

\author{Yichuan Li$^{1}$, Yonghai Gong$^{2}$, Nikolaos M. Freris$^{2}$, Petros Voulgaris$^{3}$ and Du\v{s}an Stipanovi\'{c}$^{1}$
\thanks{*This work was supported by the Ministry of Science and Technology of China under grant 2019YFB2102200, the Anhui Dept. of Science and Technology under grant 201903a05020049, the Tencent Holdings Ltd. under grant FR202003. }
\thanks{$^{1}$Coordinated Science Laboratory, University of Illinois at Urbana-Champaign, IL 61820, USA. {\tt\small yli129@illinois.edu}, {\tt\small dusan@illinois.edu}. }  %
\thanks{$^{2}$ School of Computer Science, University of Science and Technology of China, Hefei, Anhui, 230027, China. {\tt\small gongyh@mail.ustc.edu.cn}, {\tt\small nfr@ustc.edu.cn.}}
\thanks{$^{3}$Department of Mechanical Engineering, University of Nevada, Reno, NV 89557, USA. {\tt\small pvoulgaris@unr.edu}.}
}

\begin{document}

\setlength{\belowdisplayskip}{2pt} 
\setlength{\abovedisplayskip}{2pt} 
\setlength{\textfloatsep}{0pt}

\maketitle
\thispagestyle{empty}
\pagestyle{empty}

\begin{abstract}
We consider a class of distributed optimization problem where the objective function consists of a sum of strongly convex and smooth functions and a (possibly nonsmooth) convex regularizer. A multi-agent network is assumed, where each agent holds a private cost function and cooperates with its neighbors to compute the optimum of the aggregate objective. We propose a quasi-Newton Alternating Direction Method of Multipliers (ADMM) where the primal update is solved inexactly with approximated curvature information. By introducing an intermediate consensus variable, we achieve a block diagonal Hessian which eliminates the need for inner communication loops within the network when computing the update direction. We establish global linear convergence to the optimal primal-dual solution without the need for backtracking line search, under the assumption that component cost functions are strongly convex with Lipschitz continuous gradients. Numerical simulations on real datasets demonstrate the advantages of the proposed method over state of the art.
\end{abstract}

\section{INTRODUCTION} 
Distributed optimization acts as the computational engine for a wide range of applications in modern technology, including distributed control \cite{c1}, power grid management \cite{c2}, distributed machine learning \cite{c3}, and resource allocation in sensor networks \cite{c4}. We consider a classical problem of minimizing a sum of cost functions and a convex regularizer, aiming at promoting sparsity in the solution structure. In specific:
\begin{gather}
    \underset{\hat{x} \in \mathbb{R}^d}{\mathrm{minimize}}\left\{\sum_{i=1}^m f^i(\hat{x})+g(\hat{x})\right\}, \label{prob1}
\end{gather}
where $f^i:\mathbb{R}^d\to \mathbb{R}$ captures the $i$-th agent's objective of interest and $g:\mathbb{R}^d\to \mathbb{R}$ is a possibly nonsmooth regularizer, such as the weighted $\ell_1-$norm. In a network of agents, each agent aims to minimize its local cost function $f^i(\cdot)$ while communicating with its neighbors to cooperatively find the solution of the global problem. Distributed solutions are often pursued to fully utilize the computational power of the agents and reduce the amount of message passing in the network. Efficient communication is deemed especially desirable in scenarios where the number of participants $m$, and the dimension of the decision variable $d$ is large, such as in emerging applications of Cyber-Physical systems (CPS) \cite{c5}. An archetypal framework of distributed optimization introduces local decision variables at each agent, which are updated locally using private data and exchange variables with neighboring agents. A network-wide solution is achieved by means of asymptotic elimination of the consensus error, i.e., the difference between the values of neighboring agents converges to zero. Such a framework is termed distributed consensus optimization and offers significant flexibility for each agent to select appropriate updating schemes that suit its local hardware environment. 

First order methods \cite{c7}-\nocite{c8,c9}\cite{c10} constitute popular choices for solving (\ref{prob1}) in a distributed fashion. The presence of the nonsmooth regularizer prohibits gradient descent from being implemented directly as $g(\cdot)$ is not differentiable. Subgradient methods \cite{c7} invoke the notion of subdifferential set to compute descent directions while letting agents exchange information over a time-varying topology. Proximal gradient \cite{c9} accommodates the nonsmoothness by splitting the objective function and evaluating the proximal mapping associated with the nondifferentiable part. Various acceleration schemes \cite{c8}, \cite{c10} exist for first order methods, that involve storing past gradient and iterate information so as to form a momentum correction term. However, as the size of data grow the condition number of problem (\ref{prob1}) tends to get large which causes first order methods to exhibit extremely slow convergence rate. A natural remedy for the aforementioned issues is to consider second order methods \cite{c12}-\nocite{c13}\cite{c15}. Through use of function curvature information, second order methods compute  descent directions in the objective function level sets that are effective at accelerating the convergence over first order methods. Proximal Newton \cite{c15} can be considered as the second order counterpart of the proximal gradient but the associated proximal mapping increases computational burden drastically and renders distributed implementation infeasible as global information is required to compute the Newton step. Moreover, second order methods require solving a linear system to compute the Newton step which induces a computation cost of $\mathcal{O}(d^3)$. Quasi-Newton methods \cite{c16} circumvent this procedure by using a finite difference of gradients to approximate the Hessian. The reduced computation costs along with competitive performance have rendered quasi-Newton methods a desirable alternative to second order methods, especially for Large-scale optimization problems. 

Primal-dual methods \cite{c17}-\cite{onADMM} provide a different perspective for problem (\ref{prob1}) within the framework of constrained optimization, where the consensus constraint is explicitly enforced. The Alternating Direction Method of Multipliers (ADMM) \cite{onADMM} falls into this catergory where the smooth and nonsmooth parts of the objective are considered separately and iterates are computed sequentially. However, at each iteration of the ADMM, a sub-optimization problem must be solved, which typically induces an expensive computational burden. 

\textbf{Contributions}: (i) We propose BFGS-ADMM for convex composite optimization where the primal update (typically computational expensive) is replaced with a quasi-Newton step \emph{that does not require} solving a linear system of equations. Moreover, by introducing an intermediate consensus variable, we eliminate the need for inner loops within the network when computing the update direction. (ii) We establish global linear convergence rate for the proposed algorithm \emph{without backtracking}, under the assumption that private cost functions are strongly convex with Lipschitz continuous gradients. (iii) The advantage of BFGS-ADMM over first order methods is analytically established and experimentally demonstrated with numerical simulations on real datasets.

\textbf{Notation:} We denote column vectors $x\in \mathbb{R}^d$ with lower case letters and matrices $A\in\mathbb{R}^{n\times m}$ with capital letters. Superscripts are used to denote partitioned vector components and subscripts are used to denote iteration steps, e.g., $x^i_t$ denotes the value of subvector $x^i$ at step $t$. Matrix entries are denoted as $[A]^{ij}$. Unless specified otherwise, $\norm{x}$ and $\norm{A}$ denote the Euclidean norm of a vector and the corresponding induced norm of a matrix, respectively. We define the norm of a vector associated with a positive definite matrix $P\succ 0$ as $\norm{x}_P:=\sqrt{x^\top P x}$, and the set $\{1,\dots,m\}$ is abbreviated as $[m]$. The proximal mapping associated with the function $g(\cdot): \mathbb{R}^d\to \mathbb{R}$ is defined as $\textbf{prox}_{\mu g}(v) = \underset{\theta}{\mathrm{argmin}}\left\{g(\theta)+\tfrac{1}{2\mu}\norm{\theta-v}^2\right\}$.  

\section{PRELIMINARIES}
\subsection{Problem formulation}
Consider a connected network of $m$ agents captured by a $m$-th order undirected graph: $\mathcal{G}=(\mathcal{V},\mathcal{E})$, where $\mathcal{V}=\{1,\dots,m\}$ is the vertex set, $\mathcal{E}$ contains pairs $(i,j)$ if and only if agent $i$ communicates with agent $j$, and $\abs{\mathcal{E}}=n$ captures the number of communication links in the network. We denote the neighborhood of agent $i$ as $\mathcal{N}_i=\{j:(i,j)\in \mathcal{E}\}$. We further define the source and the destination matrix (arbitrary ordering suffices since both are given the same value) $\hat{A}_s,\hat{A}_d\in \mathbb{R}^{n\times m}$ respectively as follows. Each row of $\hat{A}_s$ and $\hat{A}_d$ corresponds to a communication link $\mathcal{E}_k=(i,j)$, for some $k\in[n]$, and $[\hat{A}_s]^{ki}=[\hat{A}_d]^{kj}=1$ with all other entries equal to $0$. We reformulate (\ref{prob1}) to the consensus framework by introducing local decision variables $x^i$ held by the corresponding agent $i$, along with a set of intermediate consensus variables $\{z^{ij}\}$. Specifically,
\begin{gather}
        \underset{x^i,\,\theta,\,z\,\in \mathbb{R}^d}{\mathrm{minimize}}\left\{\sum_{i=1}^m f^i(x^i)+g(\theta)\right\}\nonumber\\
    \text{s.t.} \quad 
    x^i=z^{ij}=x^j\quad \forall\, i\,\,\text{and}\, \, j\in \mathcal{N}_i,\nonumber\\
    x^l = \theta \quad \text{for some }l\in [m],\label{prob2}
\end{gather}
where the $l$-th agent is chosen arbitrarily to enforce the additional constraint for the nonsmooth regularizer decision variable. The constraints enforce consensus among the network if the underlying graph is connected and, therefore, (\ref{prob2}) is equivalent to (\ref{prob1}) in the sense that they have the same set of optimal solutions, i.e., $\hat{x}_\star=x^1_\star=\dots=x^m_\star=\theta_\star=z_\star$. Consensus can be achieved by directly enforcing $x^i=x^j$ but having a intermediate consensus variable $z^{ij}$ is \emph{crucial} in terms of decoupling the primal variable $x^i$ so that computing quasi-Newton steps does not require additional communication within the network. We provide further discussion on this design choice in Section \ref{sec3} {\it Remark 1}. We may compactly express (\ref{prob2}) using the aforementioned source and destination matrices as follows:
\begin{gather}
    \underset{x\in \mathbb{R}^{md},\theta\in \mathbb{R}^d,z\in\mathbb{R}^{nd}}{\mathrm{minimize}} \left\{F(x)+g(\theta)\right\} \nonumber\\
    \text{s.t.}\quad Ax=\begin{bmatrix}\hat{A}_s \otimes I_d \nonumber\\
    \hat{A}_d\otimes I_d\end{bmatrix}x = \begin{bmatrix}I_{nd}\nonumber\\I_{nd} \end{bmatrix}z=Dz, \nonumber\\
    Cx =\theta, \label{prob3}
\end{gather}
where $\otimes$ denotes the Kronecker product and $I$ denotes the identity matrix of the corresponding dimension. We stack local decision variables as $x\in \mathbb{R}^{md}:=[(x^1)^\top,\dots,(x^m)^\top]^\top$ and similarly for the consensus variable $z\in\mathbb{R}^{nd}:=[(z^1)^\top,\dots,(z^{n})^\top]^\top$, and matrix $A$ and $D$ is obtained by stacking matrices as shown in (\ref{prob3}). The matrix $C\in \mathbb{R}^{d\times md}:= (c^l)^\top\otimes I_d$ enforces the last line of constraint of (\ref{prob2}) using the coordinate selection vector $c^l\in \mathbb{R}^m$, whose entries are zero except the $l$-th one being 1. We present some equalities that link our constraint matrices to the topology of the graph in the following. 
\begin{gather}
    \hat{E}_s= \hat{A}_s-\hat{A}_d \quad  \hat{E}_u =\hat{A}_s+\hat{A}_d \label{incidence}\\
    \hat{E}_s^\top \hat{E}_s = \hat{L}_s \quad \hat{E}_u^\top \hat{E}_u = \hat{L}_u \label{laplacian}\\
    \hat{\Delta} = \tfrac{1}{2}(\hat{L}_s+\hat{L}_u)=\hat{A}_s^\top \hat{A}_s+\hat{A}_d^\top \hat{A}_d  \label{degree}
\end{gather}
Matrices $\hat{E}_s, \hat{E}_u\in \mathbb{R}^{n\times m}$ stand for the signed and unsigned incidence matrix of the graph, respectively, while $\hat{L}_s,\hat{L}_u\in \mathbb{R}^{m\times m}$ denotes the signed and unsigned Laplacian matrix of the graph, respectively. The diagonal degree matrix of the graph is denoted as $\hat{\Delta}$, whose entries are $\hat{\Delta}^{ii}=|\mathcal{N}_i|$.

{\it Assumption 1}. Local cost functions $f^i(\cdot)$ are twice differentiable with uniformly bounded Hessian as follows: 
\begin{gather}
    m_fI\preceq \nabla^2  f^i(\cdot) \preceq M_fI,\,\,\forall\,\,i\in[m],\label{assumption1}
\end{gather}
where $0< m_f \leq M_f<\infty$. Since $F(x)=\sum_{i=1}^m f^i(x^i)$, $\nabla^2 F(x)$ is block diagonal with the $i$-th block being $\nabla^2 f^i(x^i)$. Therefore, the same bounds apply to the Hessian of $F(x)$ as well, i.e., $m_fI\preceq \nabla^2 F(x) \preceq M_fI$. 
\subsection{Introduction to BFGS}
Quasi-Newton methods \cite{c20}\nocite{c21}-\cite{c22} constitute a popular class of methods for accelerating the convergence of optimization methods without directly invoking the Hessian. Quasi-Newton methods seek to approximate the curvature information using consecutive gradient evaluations and iterates differences. More specifically, we define the descent direction $h_t$ as
\begin{gather}
    h_t:= -B_t^{-1}\nabla f_t,
\end{gather}
where $B_t$ is a positive definite matrix that approximates the Hessian $\nabla^2 f_t$. Various schemes exist for designing algorithms that iteratively update $B_t$ including those by Davidon, Fletcher, and Powell (DFP) \cite{c20} and Broyden, Fletcher, Goldfarb, and Shannon (BFGS) \cite{c21}. \\
\indent In this paper, we focus on the BFGS scheme to select $B_t^{-1}$ as it is observed to work best in the practice both in terms of convergence speed and self-correcting capabilities \cite{c22}. The BFGS approximates the Hessian using finite differences of consecutive iterates and gradient evaluations. In specific, define the iterates difference and the gradient difference as follows:
\begin{gather*}
   s_t = x_{t+1}-x_t, \quad  d_t = \nabla f(x_{t+1})-\nabla f(x_t).
\end{gather*}
Quasi-Newton methods select $B_{t+1}$ so that the secant condition is satisfied, i.e., $B_{t+1}s_t=d_t$, which is motivated by the fact that the real Hessian satisfies this condition as $x_{t+1}$ tends to $x_t$. To ensure $B_{t+1}>0$, it must hold that $s_t^\top d_t>0$ as can seen by premultiplying the secant equation with $s_t^\top$. For convex objectives, this condition is satisfied automatically. Note that, however, the secant condition alone is not enough to specify $B_{t+1}$. To uniquely determine $B_{t+1}$, we further require its inverse to be close to the previous value $B_t^{-1}$ in the following sense: 
\begin{gather}
    \underset{B^{-1}}{\mathrm{minimize}}\quad \norm{B^{-1}-B_t^{-1}}_{\mathbf{W}}\nonumber\\
    \text{s.t.} \quad B^{-1}=(B^{-1})^\top, \quad B^{-1}d_t=s_t, \label{bfgs opti}
\end{gather}
where $\norm{A}_{\mathbf{W}}:=\norm{W^{\tfrac{1}{2}}AW^{\tfrac{1}{2}}}_{\mathbf{F}}$ stands for the weighted Frobenius norm associated with matrix $W$ whose inverse is the average Hessian \cite{c22}. A closed form solution for (\ref{bfgs opti}) can be obtained as:
\begin{gather}
    B_{t+1}^{-1} = \left(I-\tfrac{s_td_t^\top}{d_t^\top s_t}\right)B_t^{-1}\left(I-\tfrac{d_ts_t^\top}{d_t^\top s_t}\right)+\tfrac{s_ts_t^\top}{d_t^\top s_t}\label{bfgs inverse}.
\end{gather}
\subsection{Alternating Direction Method of Multipliers}
We proceed to introduce the Augmented Lagrangian associated with problem (\ref{prob2}) as follows:
\begin{gather*}
    \mathcal{L} (x,z,\theta;y,\lambda) = F(x)+g(\theta)+y^\top (Ax-Dz)\\
    +\lambda^\top (Cx-\theta)+
    \tfrac{1}{2\mu_1}\norm{Ax-Dz}^2+\tfrac{1}{2\mu_2}\norm{Cx-\theta}^2,
\end{gather*}
where $y\in \mathbb{R}^{2nd},\lambda \in \mathbb{R}^d$ are dual variables corresponding to the two linear constraints in (\ref{prob3}). The Alternating Direction Method of Multipliers (ADMM) solves (\ref{prob2}) by sequentially minimizing the Augmented Lagrangian over primal/dual variables as:
\begin{subequations}
\begin{align}
    x_{t+1} &= \underset{x}{\mathrm{argmin}}\,\,\mathcal{L}(x,z_t,\theta_t;y_t),\label{4a}\\
    z_{t+1} &= \underset{z}{\mathrm{argmin}}\,\,\mathcal{L}(x_{t+1},z,\theta_{t+1};y_t),\label{4b}\\
    \theta_{t+1} &= \underset{\theta}{\mathrm{argmin}}\,\,\mathcal{L}(x_{t+1},z_{t+1},\theta;y_t),\label{4c}\\
    y_{t+1} &= y_t+\tfrac{1}{\mu_1}(Ax_{t+1}-Dz_{t+1}) \label{4d}\\
    \lambda_{t+1} &= \lambda_t+\tfrac{1}{\mu_2}(Cx_{t+1}-\theta_{t+1}) \label{4e}
\end{align} \label{admm}%
\end{subequations}Note that (\ref{admm}) is a 3-block ADMM which, in general, is not guaranteed to converge for arbitrary values of $\mu_1,\mu_2>0$ \cite{3admm}. 
We separate the dual variables into $y$ and $\lambda$ for two reasons: (i) Since dual variables accumulate the consensus error within the network as seen in (\ref{4d}) and (\ref{4e}), and from that the fact $y$ corresponds to $2n$ constraints while $\lambda$ only involves a single constraint, it is beneficial to separate the two and design different stepsizes $\tfrac{1}{\mu_1}$ and $\tfrac{1}{\mu_2}$, which allows for extra freedom for hyperparameters and results in better performance in practice. (ii) In section \ref{sec3}, we further show that with appropriate initialization, the storage requirement for $y$ can be reduced by half. \\
\indent In each iteration of ADMM, a sub-optimization problem has to be solved in (\ref{4a}) to obtain $x_{t+1}$. This may be computationally expensive and place heavy burden on the system. We therefore propose to perform inexact optimization with the aid of an approximated Hessian using BFGS, aiming to reduce computational costs while maintaining fast convergence speed. 
\section{Algorithmic description and implementation}
We build a local model of the Augmented Lagrangian at step $t$ with respect to the primal variable $x$ as: 
\begin{gather}
    \widehat{\mathcal{L}}_t(x)= \mathcal{L}_t+\nabla_x \mathcal{L}_t^\top (x-x_t)+\tfrac{1}{2}\norm{x-x_t}^2_{H_t}, \label{model}
\end{gather}\label{sec3}where $\mathcal{L}_t:= \mathcal{L} (x_t,z_t,\theta_t;y_t,\lambda_t)$. Various algorithms can be designed by choosing different $H_t$. In this paper, we opt to use regularized BFGS approximation of the Augmented Lagrangian Hessian as follows:
\begin{gather}
    H_t = B_t+\tfrac{1}{\mu_1}\Delta+\tfrac{1}{\mu_2}C^\top C+\epsilon I_{md},\label{approx Hessian}
\end{gather}
where $\epsilon>0$ and $\Delta=A^\top A$ is the Kronecker product of the graph degree matrix $\hat{\Delta}$ and $I_d$. Note that by setting $B_t=\nabla^2 F(x_t)$, we recover the Hessian of the Augmented Lagrangian with respect to $x$, i.e., $H_t = \nabla_{xx}^2 \mathcal{L}$ in such case. The quadratic form of (\ref{model}) admits a closed form solution when minimizing over $x$, i.e., 
\begin{gather}
    x_{t+1} =x_t-H_t^{-1}\nabla_x \mathcal{L}_t \label{x}.
\end{gather}
Step (\ref{4b}) reduces to solving the following linear system of equations: 
\begin{gather}
    D^\top y_t+\tfrac{1}{\mu_1}D^\top(Ax_{t+1}-Dz_{t+1})=0. \label{z}
\end{gather}
Moreover, by completion of squares, step (\ref{4c}) is equivalent to evaluating the following proximal mapping associated with $g(\cdot)$ with parameter $\mu_2$:
\begin{gather}
    \theta_{t+1}= \textbf{prox}_{\mu_2 g}(Cx_{t+1}+\mu_2 \lambda_t). \label{theta}
\end{gather}
Dual variables are updated in verbatim as in (\ref{4d}) and (\ref{4e}). We proceed to present Lemma 1 which states that with an appropriate initialization, the intermediate variable $\{z^{ij}\}$ evolves on the manifold defined by the column space of $E_u$. Matrices without the hat denotes the Kronecker product of the corresponding matrix with the identity matrix, e.g., $E_u=\hat{E}_u\otimes I_d$ where $\hat{E}_u$ is defined in (\ref{incidence}).

\noindent\textbf{Lemma 1}. Recall the signed and unsigned incidence and Laplacian matrices in (\ref{incidence}) and (\ref{laplacian}), respectively. For zero initialization of $(x_t,y_t)$, we can decompose $y_t^\top = \begin{bmatrix}\alpha_t^\top & -\alpha_t^\top \end{bmatrix}$ for all $t\geq 0$. Moreover, the update rules (\ref{x})-(\ref{theta}), (\ref{4d}), and (\ref{4e}) can be simplified as:
\begin{subequations}
\begin{align}
        x_{t+1}&=x_t - H_t^{-1}\big\{\nabla F(x_t)+E_s^\top \alpha_t+C^\top \lambda_t+\tfrac{1}{2\mu_1}L_sx_t\nonumber  \\&+\tfrac{1}{\mu_2}C^\top (x_t^l-\theta_t)\big\},\label{14a}\\
        \theta_{t+1} &= \textbf{prox}_{\mu_2 g}(x_{t+1}^l+\mu_2 \lambda_t),\label{14b}\\
        \alpha_{t+1} &= \alpha_t+\tfrac{1}{2\mu_1} E_s x_{t+1},\label{14c} \\
        \lambda_{t+1} &= \lambda_t+\tfrac{1}{\mu_2}(x_{t+1}^l-\theta_{t+1}).\label{14d}
\end{align}\label{updates}
\end{subequations}
{\it Proof} : See Appendix A. \\
Once $\theta_{t+1}$ is obtained, updates (\ref{14c}) and (\ref{14d}) can be performed in parallel and we have effectively transformed a 3-block ADMM in (\ref{admm}) to a 2-block ADMM, which converges under more general settings \cite{c23}. \\
\indent {\it Remark 1}: The reason for introducing the consensus variable $z$ in (\ref{prob2}) is to decouple $x^i$'s so that the Hessian of the Augmented Lagrangian has a block diagonal structure, cf.(\ref{approx Hessian}), which is instrumental for computing the quasi-Newton step {\it without additional communication} within the network once the gradient of the Augmented Lagrangian is obtained. Note that there is no need to invert any matrix or solve a linear system, and the inverse in (\ref{updates}), is only notational. Instead, we directly approximate each block $(B_{t+1}^{-1})^{ii}$, corresponding to the agent $i$ and form $(H_{t+1}^{-1})^{ii}$ as explained in the following. Since $B_{t+1}^{-1}$ is block diagonal, we approximate $(B_{t+1}^{-1})^{ii}$ using (\ref{bfgs inverse}) with $d_t=\nabla f^i(x^i_{t+1})-\nabla f^i(x^i_t)$ and $s_t=x_{t+1}^i-x_t^i$. Moreover, $H_{t+1}-B_{t+1}$ is diagonal and constant (time invariant) from (\ref{approx Hessian}), which allows for eigendecomposition with {\it components only depending on local information}. Specifically, using the Sherman-Morrison formula and defining $C^i_1=(B_{t+1}^{-1})^{ii}$, we obtain an iterative formula for $(H_{t+1}^{-1})^{ii}=C^i_{d+1}$ as follows:
\begin{gather}
C^i_{k+1}=C^i_k-\tfrac{C^i_kv_k^{i}(v^i_k)^\top C^i_k}{1+(v^k)^\top C^i_kv^k},  \label{H update}
\end{gather}
with constant vector $v^i_k$ specified as:
\begin{gather}
    v^i_k = \left(\tfrac{\abs{\mathcal{N}_i}}{\mu_1}+\tfrac{\mathbf{1}_{il}}{\mu_2}+\epsilon\right)e_k ,
\end{gather}
where $\mathbf{1}_{il}=1$ if $i=l$ (i.e., if it is the $l$-th agent who performs the proximal mapping associated with the nonsmooth regularizer $g(\cdot)$) and is zero otherwise, and $e_k\in\mathbb{R}^d$ is a constant vector with all entries equal to $0$ except for the $k$-th one. If consensus is enforced directly as $x^i=x^j$, the resulting Hessian of the Augmented Lagrangian would involve the graph Laplacian matrix as in \cite{c13} which couples $x^i$ with its neighbors. In other words, computing quasi-Newton steps would induce multiple inner communication rounds within the network at each iteration, to approximate the descent direction to a desired accuracy.\\
\indent We now present a first order variant of BFGS-ADMM, which only uses first order information, so as to demonstrate the merits of using the BFGS approximation, both theoretically and experimentally. Note that the only part of $H_t$ in (\ref{approx Hessian}) that depends on the iterate counter is $B_t$, which is an approximation of $\nabla^2 F(x_t)$. Therefore, a first order approximation can be derived by setting $B_t=0$ in which case $H_t$ in (\ref{approx Hessian}) is a constant diagonal matrix. In other words, agent $i$ selects a step size equal to:
\begin{gather}
    s_i = \left(\tfrac{\abs{\mathcal{N}_i}}{\mu_1}+\tfrac{\mathbf{1}_{il}}{\mu_2}+\epsilon\right)^{-1}. \label{first order}
\end{gather}
Therefore the first order variant updates the primal vector $x_{t+1}$ as:
\begin{gather}
    x_{t+1}=x_t -\textbf{diag}[s_i]\big\{\nabla F(x_t)+E_s^\top \alpha_t+C^\top \lambda_t+\tfrac{1}{2\mu_1}L_sx_t\nonumber  \\
    +\tfrac{1}{\mu_2}C^\top (x_t^l-\theta_t)\big\}, \label{1admm}
\end{gather}
\begin{algorithm}[t]
    \caption{BFGS-ADMM update} 
    Zero initialization for all variables. Hyperparameters $\mu_1,\mu_2$, and $\epsilon$. Initialization for $B_0^{-1}=aI_{md}$ for some $a>0$.
    \begin{algorithmic}[1]
    \For{$t=0,1,2,\ldots$}
    \For{$i\in[m]$}
    \State Compute $h_t^i=\nabla f^i(x_t^i)+\tfrac{1}{2\mu_1}\sum_{j\in \mathcal{N}^i} (x_t^i-x_t^j) + \phi_t^i$.
            \If{$i=l$}
                \State Update $h_t^i\leftarrow h_t^i +\tfrac{1}{\mu_2}(x_t^i-\theta_{t})+\lambda_t$.
            \EndIf
    \State Compute $u^i_t=(H_t^{-1})^{ii}h_t^i$.
    \State \textbf{Primal update:} $x_{t+1}^i=x_t^i-u_t^i$.
    \State \textbf{Dual update:} $  \phi_{t+1}^i = \phi_t^i+\tfrac{1}{2\mu_1}\sum_{j\in \mathcal{N}^i}(x_{t+1}^i-x_{t+1}^j)$.
    \If{$i=l$}
    \State Compute $\theta_{t+1}=\textbf{prox}_{\mu_2 g}(x^i_{t+1}+\mu_2 \lambda_t)$.
    \State Compute $\lambda_{t+1} = \lambda_t +\tfrac{1}{\mu_2}(x^i_{t+1}-\theta_{t+1})$.
    \EndIf
    \State Set $s_t=-u_t^i$ and $d_t= \nabla f^i(x^i_{t+1})-\nabla f(x^i_t)$.
    \State Update $(B_{t+1}^{-1})^{ii}$ using (\ref{bfgs inverse}).
    \State Update $(H_{t+1}^{-1})^{ii}$ using (\ref{H update}).
    \EndFor
    \EndFor
\end{algorithmic}
\end{algorithm}where $\textbf{diag}[s_i]\in\mathbb{R}^{md}$ is the diagonal matrix formed by $s_i I_{md}$, and the dual variables $(\theta_{t+1},\alpha_{t+1},\lambda_{t+1})$ are updated in the exact same way as in the BFGS-ADMM using (\ref{14b})-(\ref{14c}). Note that our first order variant encapsulates existing first order methods EXTRA \cite{extra} as a special case by setting $s_i=\tfrac{1}{\epsilon}$. Indeed, when $g(\cdot)=0$, the proximal mapping and the associated dual variables in (\ref{updates}) are not needed and we obtain:
$
    x_{t+1}=x_t-\tfrac{1}{\epsilon}\left[\nabla F(x_t)+E_s^\top \alpha_t+\tfrac{1}{2\mu_1}L_sx_t\right]. \label{extra1}
$
Taking the difference between the $x_{t+1}$ update and the $x_t$ update defined this way, and substituting the dual update for $\alpha_t$ using (\ref{14c}), we obtain:
\begin{gather}
    x_{t+1}=2(I-\tfrac{1}{2\mu_1}L_s)x_t-(I-\tfrac{1}{2\mu_1}L_s)x_{t-1}\nonumber\\-\tfrac{1}{\epsilon}(\nabla F(x_t)-\nabla F(x_{t-1})). \label{extra2}
\end{gather}
With appropriate choices of $\mu_1,\mu_2,\epsilon$, the update defined in (\ref{extra2}) is equivalent to EXTRA. We present the comparison between BFGS-ADMM and its first order variant in terms of approximating the exact ADMM in the next section.\\
\indent We proceed to present the distributed implementation of BFGS-ADMM in \textbf{Algorithm 1}. Note that in the primal update (\ref{14a}), dual variables are invoked in the form of $E_s^\top \alpha_t$. By defining $\phi_t= E_s^\top \alpha_t\in \mathbb{R}^{md}$ and pre-multiplying both sides of (\ref{14c}) with $E_s^\top$, we obtain an equivalent algorithm that allows for efficient distributed implementation. We let each agent hold the corresponding pair $(x^i,\phi^i)$ while the $l$-th agent additionally holds $(\theta,\lambda)$. At each iteration, each agent $i$ begins with calculating $h_t^i$ using $x_t^j$ from its neighbors as in step 3, where the $l$-th agent performs an additional update in step 5. Note that step 7 does not require solving a linear system since $(H_t^{-1})^{ii}$ is directly approximated using (\ref{H update}). Primal and dual updates are performed at each agent as in step 8 and step 9, respectively. The $l$-th agent evaluates an additional proximal mapping in step 11 to update $\theta_{t+1}$, and updates $\lambda_{t+1}$ in step 12. Finally, each agent $i$ updates its $(B_{t+1}^{-1})^{ii}$ and $(H_{t+1}^{-1})^{ii}$ using (\ref{bfgs inverse}) and (\ref{H update}), respectively.

\section{CONVERGENCE ANALYSIS}
Linear convergence rate for ADMM is well established \cite{onADMM}. Since our goal here is to reduce the computational cost of ADMM, the best one can hope for is to reduce the gap between the approximation and the standard ADMM, while maintaining the linear convergence rate. In \textbf{Lemma 4}, we demonstrate the advantage of using BFGS approximation over first order methods by showing that the aforementioned gap decreases faster for BFGS-ADMM. We first state an additional assumption on the nonsmooth regularizer $g(\cdot)$ and the KKT conditions associated with (\ref{prob3}) expressed in the primal optimal pair $(x_\star,\theta_\star)$ and the dual optimal pair $(\alpha_\star,\lambda_\star)$. Note that we have eliminated $z_\star$, as it is shown in \textbf{Lemma 1} that it is a function of $x_\star$.

{\it Assumption 2}: The function $g(\cdot):\mathbb{R}^d\to \mathbb{R}$ is proper, closed, and convex. Equivalently, $\forall\,x,y\in\mathbb{R}^d$, the following inequality holds:
\begin{gather}
    (x-y)^\top (\partial g(x)-\partial g(y)) \geq 0,
\end{gather}
where $\partial g(\cdot)$ denotes the subdifferential set. We proceed to state a non-restrictive assumption that can be easily achieved by regularization.

{\it Assumption 3}: There exists positive constant $0<\nu<\infty$ such that the eigenvalues of $B_t$ are upper bounded, i.e, 
\begin{gather*}
    B_t < \nu I.
\end{gather*}
\indent {\it Remark 2}: The upper bound described above can be achieved by setting $B_t^{-1} = \hat{B}_t^{-1}+\tfrac{1}{\nu}I$, where $\hat{B}_t$ comes from the BFGS update formula. Since $\hat{B}_t^{-1} >0$ with positive definite initialization, we have $B_t^{-1}> \tfrac{1}{\nu}I$, i.e., $B_t <\nu I$. \\
\noindent \textbf{KKT conditions}:
\begin{subequations}
\begin{align}
    \nabla F(x_\star)+E_s^\top \alpha_\star+C^\top \lambda_\star &=0 \tag{KKTa} \\
    \partial g(\theta_\star)-\lambda_\star &\ni 0 \tag{KKTb}\\
    E_s x_\star &= 0 \tag{KKTc} \\
    x_\star^l -\theta_\star & =0  \tag{KKTd}
\end{align}
\end{subequations}

We proceed to state a technical lemma that describes an inclusion relationship between the subdifferential set of $\partial g(\theta_{t+1})$ and the dual variable $\lambda_{t+1}$. 

\noindent \textbf{Lemma 2}. Consider the update $\theta_{t+1}=\textbf{prox}_{\mu_2 g}(x_{t+1}^l+\mu_2\lambda_t)$ in (\ref{14b}). It holds that:
\begin{gather}
    \lambda_{t+1} \in \partial g(\theta_{t+1}).
\end{gather}
{\it Proof} : See Appendix B.

Recall the first order variant of the BFGS-ADMM where the primal variable is updated as in (\ref{1admm}) while BFGS-ADMM updates $x_{t+1}$ as in (\ref{14a}), and both update dual variables as in (\ref{14b})-(\ref{14c}). We proceed to present results which capture their similarities and differences in the following.

\noindent\textbf{Lemma 3}. If {\it Assumptions 1-2} hold, then the iterates generated by BFGS-ADMM and its first order variant both satisfy the following equation:
\begin{gather}
    e_t+\nabla F(x_{t+1})-\nabla F(x_\star)+\tfrac{1}{2\mu_1}(L_u+\epsilon I)(x_{t+1}-x_t)\nonumber\\
    +C^\top \big\{\lambda_{t+1}-\lambda_\star+\tfrac{1}{\mu_2}(\theta_{t+1}-\theta_t) \big\}+E_s^\top (\alpha_{t+1}-\alpha_\star)=0, \label{lemma3 eqn}
\end{gather}
where for BFGS-ADMM,
\begin{gather}
e_t=\nabla F(x_t)+B_t(x_{t+1}-x_t)-\nabla F(x_{t+1}).\label{error}
\end{gather}
for its first order variant,
\begin{gather}
    e_t = \nabla F(x_t)-\nabla F(x_{t+1}).  \label{error2}
\end{gather}
{\it Proof} : See Appendix C.\\
By comparing (\ref{error}) and (\ref{error2}), we see that the difference between using BFGS vs. first order approximation lies in how $\nabla F(x_{t+1})$ is approximated at step $t$. Moreover, if the sub-optimization problem is solved exactly as in (\ref{4a}), by optimality condition, equation (\ref{lemma3 eqn}) holds with $e_t=0$.\\
\noindent \textbf{Lemma 4}. Consider the bound for $\nabla^2 F(x)$ in (\ref{assumption1}) and the $e_t$ introduced in \textbf{Lemma 3}. If {\it Assumption 3} holds,  an upper bound can be established as follows. \\
For BFGS-ADMM,
\begin{align}
    \norm{e_t}&\leq \norm{B_{t+1}-B_t}\norm{x_{t+1}-x_t}  \label{bfgs error1}
\end{align}
For its first order variant as in (\ref{1admm}), 
\begin{gather}
    \norm{e_t}\leq M_f\norm{x_{t+1}-x_t}. \label{1st order error}
\end{gather}
{\it Proof} : See Appendix D.\\
\indent {\it Remark 3}: \textbf{Lemma 4} is a standalone result whose proof does not require the convergence of the algorithm. In \textbf{Theorem 1}, we establish that BFGS-ADMM converges linearly which means $\lim_{t\to \infty}\norm{B_{t+1}-B_t}=0$, since both $d_t,s_t$ approach $0$ as $t\to \infty$ in (\ref{bfgs inverse}). Therefore, $\norm{e_t}=o(\norm{x_{t+1}-x_t})$ for the BFGS-ADMM while $\norm{e_t}=O(\norm{x_{t+1}-x_t})$ for its first order variant. The variable $e_t$ captures the approximation error vs. exact solution of the primal problem ($e_t=0$ for this case), whence \textbf{Lemma 4} unravels that the proposed has superior tracking of the exact ADMM, which is the key attribute to yield superior linear convergence rate over its first-order variant. \\
\noindent \textbf{Lemma 5}. Consider the update in (\ref{14c}) and (\ref{14d}). With zero initialization, the stacked vector $[\alpha^\top_t\,\, \lambda_t^\top]^\top \in \mathbb{R}^{(n+1)d}$ lies in the column space of $X^\top=[E_s^\top \,\, C^\top]^\top \in \mathbb{R}^{(n+1)d\times md}$. Moreover, there exists a unique optimal $[\alpha_\star^\top\,\, \lambda_\star^\top]^\top$ in the column space of $X^\top$ and, denoting $\sigma^+_{\mathrm{min}}$ as the smallest positive eigenvalue of $X^\top X$, the following inequality holds: 
\begin{gather}
\sigma^+_{\mathrm{min}} \left(\norm{\alpha_{t+1}-\alpha_\star}^2+\norm{\lambda_{t+1}-\lambda_\star}^2 \right) \leq \nonumber \\
    \norm{ E_s^\top (\alpha_{t+1}-\alpha_\star)+C^\top (\lambda_{t+1}-\lambda_\star)}^2. \label{lemma5}
\end{gather}
{\it Proof} : See Appendix E. \\
\indent We introduce the following scaling matrix defined in terms of the hyperparameters $\epsilon, \mu_1, \mu_2$ and the graph topology (as captured by $L_u$), 
$
    \mathcal{H} = \begin{bmatrix}
    \tfrac{1}{2\mu_1}L_u+\epsilon I & 0 & 0& 0 \\
    0 & \tfrac{1}{\mu_2} & 0 &0 \\
    0 & 0 & 2\mu_1 &0 \\
    0 & 0& 0 & \mu_2
    \end{bmatrix}.
$
\noindent \textbf{Theorem 1}. Consider the update in (\ref{updates}). Define $u^\top = [x^\top,\theta^\top,\alpha^\top,\lambda^\top]$ and $u_\star$ the corresponding optimum. Denote the smallest positive eigenvalue of $\begin{bmatrix}
E_s\\ C
\end{bmatrix}\begin{bmatrix}
E_s^\top C^\top
\end{bmatrix}$ as $\sigma^{+}_{\mathrm{min}}$, the smallest and largest eigenvalue of $L_u$ as $\sigma^G_{\mathrm{min}}$ and $\sigma^G_{\mathrm{max}}$, respectively. Denoting $\kappa=\tfrac{M_f}{m_f}$, for arbitrary constants $\beta,\gamma,\phi,\rho>1$, and $\zeta \in (\tfrac{m_f+M_f}{2m_fM_f},\tfrac{\epsilon}{4M_f^2})$, $\epsilon>2(m_f+M_f)\kappa$. the iterates converge linearly as follows:
\begin{gather}
    \norm{u_{t+1}-u_\star}^2_\mathcal{H} \leq \tfrac{1}{1+\delta}\norm{u_t-u_\star}^2_{\mathcal{H}} \label{theorem1},
\end{gather}
where 
\begin{gather}
        \delta = \min\bigg\{\left(\tfrac{2m_fM_f}{m_f+M_f}-\tfrac{1}{\zeta}\right)\left(\tfrac{2\mu_1\mu_2}{\mu_2(\sigma_{\max}^G+2\mu_1\epsilon)+\rho\mu_1}\right), \tfrac{\rho-1}{\rho},\nonumber\\ \tfrac{(\beta-1)\sigma^+_{\mathrm{min}}}{\mu\beta}\left(\tfrac{\tfrac{\sigma_{\min}^G}{2\mu_1}+\epsilon-\zeta\tau^2}{\tfrac{\sigma_{\max}^G}{2\mu_1}+\epsilon+\tfrac{\gamma\tau^2(\beta-1)}{\sigma^+_{\mathrm{min}}(\gamma-1)}}\right), \nonumber\\
    \tfrac{2\sigma^+_{\mathrm{min}}}{(m_f+M_f)\mu\beta\psi}\left(\tfrac{1}{\gamma}\right), \tfrac{\sigma^+_{\mathrm{min}}\mu_2(\psi-1)}{\mu\beta\psi}\left(\tfrac{1}{\gamma}\right)
    \bigg\}.\label{choose}
\end{gather} 
{\it Proof:} Since $F(x)$ is strongly convex with parameter $m_f$ and the gradient $\nabla F(x)$ is Lipschitz continuous with parameter $M_f$,
\begin{gather}
    \tfrac{m_fM_f}{m_f+M_f}\norm{x_{t+1}-x_\star}^2 
    +\tfrac{1}{m_f+M_f}\norm{\nabla F(x_{t+1})-\nabla F(x_\star)}^2  \leq \nonumber\\
    (x_{t+1}-x_\star)^\top  (\nabla F(x_{t+1})-\nabla F(x_\star)) . \label{coer}
\end{gather}
Using \textbf{Lemma 3}, we rewrite the \textbf{RHS} (right hand side) of (\ref{coer}) as 
\begin{gather}
    \textbf{RHS} =-(x_{t+1}-x_\star)^\top \tfrac{1}{2\mu_1}(L_u+\epsilon I)(x_{t+1}-x_t)\nonumber\\
    -(x_{t+1}^l-x_\star^l)^\top \left[\lambda_{t+1}-\lambda_\star+\tfrac{1}{\mu_2}(\theta_{t+1}-\theta_t)\right]\nonumber \\
    -(x_{t+1}-x_\star)^\top E_s^\top (\alpha_{t+1}-\alpha_\star)-(x_{t+1}-x_\star)^\top e_t, \label{p3}
\end{gather}
where we have used the fact that $(x_{t+1}-x_\star)^\top C^\top = (x_{t+1}^l-x_\star^l)$. From the dual update (\ref{14c}) and (KKTc), it holds that 
\begin{gather}
    (x_{t+1}-x_\star)^\top E_s^\top= 2\mu_1 (\alpha_{t+1}-\alpha_t)^\top . \label{p1}
\end{gather}
Moreover, using (\ref{14d}) and (KKTd), we obtain 
\begin{gather}
    (x_{t+1}^l-x_\star^l)^\top = \mu_2(\lambda_{t+1}-\lambda_t)^\top \label{p2} +(\theta_{t+1}-\theta_\star)^\top.
\end{gather}
After substituting (\ref{p1}) and (\ref{p2}) into (\ref{p3}), we obtain
\begin{gather}
    \textbf{RHS} = -(x_{t+1}-x_\star)^\top \tfrac{1}{2\mu_1}(L_u+\epsilon I)(x_{t+1}-x_t)\nonumber \\
    -\mu_2(\lambda_{t+1}-\lambda_t)^\top (\lambda_{t+1}-\lambda_\star)-(\lambda_{t+1}-\lambda_t)^\top(\theta_{t+1}-\theta_t)\nonumber\\
    -(\theta_{t+1}-\theta_\star)^\top(\lambda_{t+1}-\lambda_\star)-\tfrac{1}{\mu_2}(\theta_{t+1}-\theta_\star)^\top(\theta_{t+1}-\theta_t)\nonumber\\
    -2\mu_1(\alpha_{t+1}-\alpha_t)^\top (\alpha_{t+1}-\alpha_\star)-(x_{t+1}-x_\star)^\top e_t.
\end{gather}
Using \textbf{Lemma 2}, we have $(\lambda_{t+1}-\lambda_t)^\top(\theta_{t+1}-\theta_t)\in (\partial g(\theta_{t+1})-\partial g(\theta_t))^\top (\theta_{t+1}-\theta_t)\geq 0$, where the inequality follows from \textbf{Lemma 2} and \textbf{Assumption 1}. Similarly, $(\theta_{t+1}-\theta_\star)^\top (\lambda_{t+1}-\lambda_\star)\geq0$. Therefore, we obtain 
\begin{gather}
        \textbf{RHS} \leq -(x_{t+1}-x_\star)^\top \tfrac{1}{2\mu_1}(L_u+\epsilon I)(x_{t+1}-x_t)\nonumber \\
    -\mu_2(\lambda_{t+1}-\lambda_t)^\top (\lambda_{t+1}-\lambda_\star)
    -(\theta_{t+1}-\theta_\star)^\top\tfrac{1}{\mu_2}(\theta_{t+1}-\theta_t)\nonumber\\
    -2\mu_1(\alpha_{t+1}-\alpha_t)^\top (\alpha_{t+1}-\alpha_\star)-(x_{t+1}-x_\star)^\top e_t.\label{p4}
\end{gather}
Note that $-2(a-b)^\top  (a-c)=\norm{b-c}^2-\norm{a-b}^2-\norm{a-c}^2$ holds true for any $(a,b,c)$. Multiplying both sides of (\ref{p4}) by $2$ and using the aforementioned identity for the inner product terms, and considering the concatenation $u^\top = [x^\top,\theta^\top,\alpha^\top,\lambda^\top]$ we obtain:
\begin{gather}
    \tfrac{2mM}{m+M}\norm{x_{t+1}-x_\star}^2+\tfrac{2}{m+M}\norm{\nabla F(x_{t+1})-\nabla F(x_\star)}^2\nonumber \\
    +\norm{x_{t+1}-x_t}^2_{\Tilde{L}_u}+\mu_2\norm{\lambda_{t+1}-\lambda_t}^2+\tfrac{1}{\mu_2}\norm{\theta_{t+1}-\theta_t}^2\nonumber\\
    +2\mu_1\norm{\alpha_{t+1}-\alpha_t}^2+2(x_{t+1}-x_\star)^\top e_t \leq\nonumber \\
    \norm{u_t-u_\star}^2_{\mathcal{H}}-\norm{u_{t+1}-u_\star}^2_{\mathcal{H}} \label{p6}
\end{gather}
where $\Tilde{L}_u= \tfrac{1}{2\mu_1}(L_u+\epsilon I)$. To establish linear convergence as in (\ref{theorem1}), it suffices to show that for some $\delta>0$, the following holds: 
$
    \delta\norm{u_{t+1}-u_\star}^2_{\mathcal{H}} \leq \norm{u_t-u_\star}^2_{\mathcal{H}}-\norm{u_{t+1}-u_\star}^2_{\mathcal{H}}
$.
Moreover, for any $\zeta>0$, the inequality holds: $2(x_{t+1}-x_\star)^\top e_t \geq -\tfrac{1}{\zeta}\norm{x_{t+1}-x_\star}^2-\zeta\norm{e_t}^2$. Therefore, it is sufficient to show 
\begin{gather}
     \tfrac{2m_fM_f}{m_f+M_f}\norm{x_{t+1}-x_\star}^2+\tfrac{2}{m_f+M_f}\norm{\nabla F(x_{t+1})-\nabla F(x_\star)}^2\nonumber \\
    +\norm{x_{t+1}-x_t}^2_{\Tilde{L}_u}+\mu_2\norm{\lambda_{t+1}-\lambda_t}^2+\tfrac{1}{\mu_2}\norm{\theta_{t+1}-\theta_t}^2\nonumber\\
    +2\mu_1\norm{\alpha_{t+1}-\alpha_t}^2-\tfrac{1}{\zeta}\norm{x_{t+1}-x_\star}^2-\zeta\norm{e_t}^2 \geq \nonumber \\
    \delta\norm{u_{t+1}-u_\star}^2_{\mathcal{H}} . \label{p12}
\end{gather}
We proceed to establish this bound. From \textbf{Lemma 3}, it holds: 
\begin{gather}
    E_s^\top (\alpha_{t+1}-\alpha_\star)+C^\top (\lambda_{t+1}-\lambda_\star) = \nonumber \\
    -\big\{\nabla F(x_{t+1})-\nabla F(x_\star)+\tfrac{1}{\mu_1}(L_u+\epsilon I)(x_{t+1}-x_t)\nonumber\\
    +\tfrac{1}{\mu_2}C^\top(\theta_{t+1}-\theta_t)+e_t \big\}. \label{p7}
\end{gather}
Note that $\norm{a}^2\leq \tfrac{\beta}{\beta-1}\norm{b}^2+\beta\norm{c}^2$ holds true for any $\beta>1$. After applying this inequality three times with arbitrary constant $\beta,\gamma,\psi>1$ for (\ref{p7}), we obtain
\begin{gather}
   \norm{ E_s^\top (\alpha_{t+1}-\alpha_\star)+C^\top (\lambda_{t+1}-\lambda_\star)}^2 \leq \nonumber \\ \tfrac{\beta}{\beta-1} \norm{x_{t+1}-x_t}^2_{\Tilde{L}_u^2} 
   +\tfrac{\beta\gamma}{\gamma-1}\norm{e_t}^2+\tfrac{\beta\gamma\psi}{(\mu_2)^2(\psi-1)}\norm{\theta_{t+1}-\theta_t}^2\nonumber \\
   +\beta \gamma \psi \norm{\nabla F(x_{t+1})-\nabla F(x_\star)}^2. \label{p8}
\end{gather}
Consider the error term in \textbf{Lemma 4}, $
     e_t=B_{t+1}-B_t=-\tfrac{B_ts_ts_t^\top B_t}{s_t^\top B_ts_t}+\tfrac{d_td_t^\top}{d_t^\top s_t}.
$
Therefore, 
$
    \norm{B_{t+1}-B_t} \leq \tfrac{\norm{d_t}^2}{d_t^\top s_t}+\tfrac{\norm{B_ts_t}^2}{s_t^\top B_ts_t}
$.
Since $F(x)$ is strongly convex with $m_f$ and the gradient is Lipschitz continuous with parameter $M_f$, the following inequality holds: 
\begin{gather}
    d_t^\top s_t \geq m_f \norm{s_t}^2,\quad \text{and} \quad
    d_t^\top s_t \geq \tfrac{1}{M_f}\norm{d_t}^2. \label{ple4a}
\end{gather}
By setting $q_t=B_t^{\tfrac{1}{2}}s_t$, we obtain
$
  \tfrac{\norm{B_ts_t}^2}{s_t^\top B_t s_t}= \tfrac{q_t^\top B_t q_t}{q_t^\top q_t}\leq \lambda_{\mathrm{max}}(B_t)<\nu. 
$
Therefore,
$
    \norm{e_t}\leq \tau \norm{x_{t+1}-x_t},
$
where $\tau = \min\left\{M,\tfrac{\norm{d_t}^2}{\norm{s_t}^2}\right\}+\nu$.
From \textbf{Lemma 5}, we have 
\begin{gather}
    \sigma^+_{\mathrm{min}} \left(\norm{\alpha_{t+1}-\alpha_\star}^2+\norm{\lambda_{t+1}-\lambda_\star}^2 \right) \leq \nonumber \\
    \norm{ E_s^\top (\alpha_{t+1}-\alpha_\star)+C^\top (\lambda_{t+1}-\lambda_\star)}^2.\label{p9}
\end{gather}
Denoting $\mu =\max\{2\mu_1,\mu_2\}$ and combining (\ref{p8})-(\ref{p9}), we have: 
\begin{gather}
    2\mu_1\norm{\alpha_{t+1}-\alpha_\star}^2+\mu_2 \norm{\lambda_{t+1}-\lambda_\star}^2 \leq \nonumber\\
    \tfrac{\mu\beta}{\sigma^+_{\mathrm{min}}(\beta-1)} \norm{x_{t+1}-x_t}^2_{\Tilde{L}_u^2} 
   +\tfrac{\mu\beta\gamma\tau^2}{\sigma^+_{\mathrm{min}}(\gamma-1)}\norm{x_{t+1}-x_t}^2\nonumber\\
   +\tfrac{\mu\beta\gamma\psi}{(\mu_2)^2\sigma^+_{\mathrm{min}}(\psi-1)}\norm{\theta_{t+1}-\theta_t}^2\nonumber\\
   +\tfrac{\mu\beta \gamma \psi}{\sigma^+_{\mathrm{min}}} \norm{\nabla F(x_{t+1})-\nabla F(x_\star)}^2. \label{p10}
\end{gather}
Furthermore, from (\ref{14d}) and (KKTd), it holds that 
$
    \theta_{t+1}-\theta_\star = -\mu_2(\lambda_{t+1}-\lambda_t)+x_{t+1}^l-x_\star^l
$.
Using the same technique as in deriving (\ref{p8}), we obtain for arbitrary $\rho>1$:
\begin{gather}
    \norm{\theta_{t+1}-\theta_\star}\leq \rho \norm{x_{t+1}^l-x_\star^l}^2+\tfrac{(\mu_2)^2\rho}{\rho-1}\norm{\lambda_{t+1}-\lambda_t}^2. \label{p11}
\end{gather}
Using (\ref{p10}) and (\ref{p11}), it suffices to show that the following inequality holds true for some $\delta>0$,
\begin{gather}
    \delta\norm{u_{t+1}-u_\star}^2_{\mathcal{H}} \leq 
     \tfrac{2m_fM_f}{m_f+M_f}\norm{x_{t+1}-x_\star}^2\nonumber\\
     +\tfrac{2}{m_f+M_f}\norm{\nabla F(x_{t+1})-\nabla F(x_\star)}^2\nonumber \\
    +\norm{x_{t+1}-x_t}^2_{\Tilde{L}_u}+\mu_2\norm{\lambda_{t+1}-\lambda_t}^2+\tfrac{1}{\mu_2}\norm{\theta_{t+1}-\theta_t}^2\nonumber\\
    +2\mu_1\norm{\alpha_{t+1}-\alpha_t}^2-\tfrac{1}{\zeta}\norm{x_{t+1}-x_\star}^2-\zeta\tau^2\norm{x_{t+1}-x_t}^2 \nonumber,
\end{gather}
which is satisfied if $\delta$ is chosen as in (\ref{choose}).\QEDB

\section{EXPERIMENTS}

\begin{figure}[t]
  \centering
  \begin{subfigure}[b]{0.49\linewidth}
    \includegraphics[width=\textwidth]{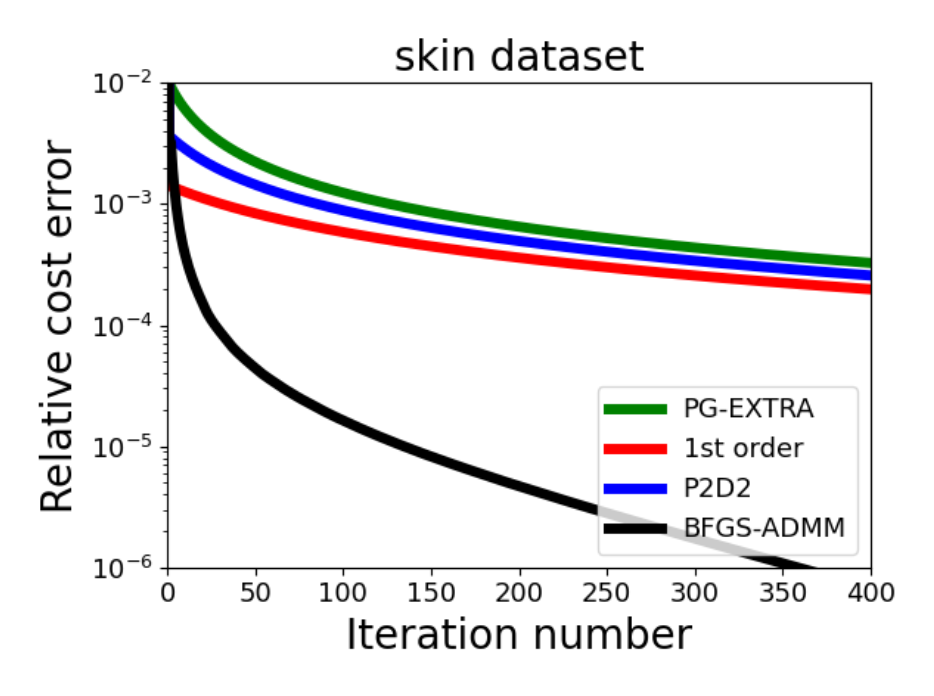}
    \caption{$m = 20,\,d=3$.}
  \end{subfigure}
  \begin{subfigure}[b]{0.49\linewidth}
    \includegraphics[width=\textwidth]{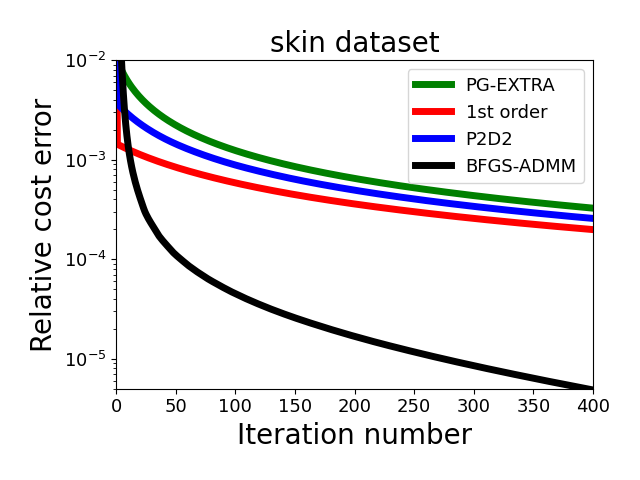}
    \caption{$m=40,\,d=3$.}
  \end{subfigure}
  \caption{skin\_noskin dataset.}
  \label{fig1}
  \begin{subfigure}[b]{0.49\linewidth}
    \includegraphics[width=\textwidth]{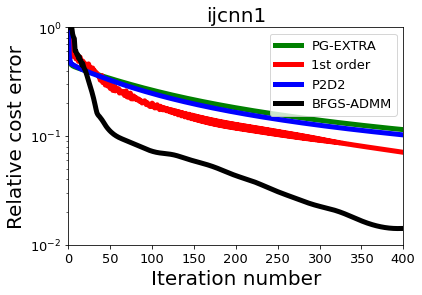}
    \caption{$m=20,\,d=22$.}
  \end{subfigure}
  \begin{subfigure}[b]{0.49\linewidth}
    \includegraphics[width=\textwidth]{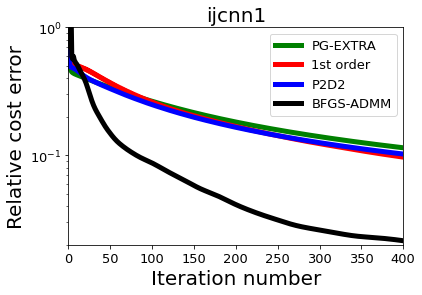}
    \caption{$m=40,\,d=22$.}
  \end{subfigure}
  \caption{ijcnn1 dataset.}
  \label{fig2}
\end{figure}

In this section we present some numerical experiments of the proposed method compared to distributed first order methods: P2D2 \cite{c24} and PG-EXTRA \cite{c8}, for solving the following distributed logistic regression problem:
\begin{gather*}
\underset{x\in \mathbb{R}^d}{\mathrm{minimize}}\,\,F(x)=\left\{\frac{1}{m}\sum_{i=1}^m f^i(x)  +\gamma\norm{x}_1  \right\},
\end{gather*}
where $
    \setlength{\belowdisplayskip}{10pt}
    f^i(x) = \tfrac{1}{m_i}\sum_{j=1}^{m_i} \left[\ln(1+e^{w_j^\top x})+(1-y_j)w_j^\top x\right],$
and $m_i$ is the number of data points held at each agent and $(w_j,y_j)\in \mathbb{R}^{d}\times \{0,1\}$ are training samples of dimension $d$ and binary labels, respectively. We consider two real datasets from LIBSVM\footnote{https://www.csie.ntu.edu.tw/~cjlin/libsvm/}: the skin\_nonskin dataset and the ijcnn1 dataset. We take 5,000 data points from each dataset with dimension $d=3$ and dimension $d=22$, respectively. A random binomial graph with $m$ agents is drawn (each edge is drawn independently from a Bernoulli($p$) distribution; $p$=0.2 was used in all cases). The mixing matrix for P2D2 is generated using the Metropolis rule while the mixing matrix of PG-EXTRA is generated using the Laplacian based constant edge weight matrix. In all cases, the $\ell_1-$norm weight $\gamma=2\times 10^{-6}$. We plot the average relative cost error defined as 
$
    \tfrac{\tfrac{1}{m}\sum_{i=1}^m F(x_t^i)-F(x_\star)}{\tfrac{1}{m}\sum_{i=1}^m F(x_0^i)-F(x_\star)}
$ for different network sizes. As it can be seen in Figure 1 and 2, the BFGS-ADMM method demonstrates significant speed up in both datasets compared to other methods. Note that our first order variant also shows advantages compared to other first order methods, due to the fact that in our first order approximation, step size selection also takes into account the number of neighbors each agent has as in (\ref{first order}). 

\bibliography{BFGS_full.bib}
\bibliographystyle{IEEEtran}

\appendix
\noindent \textbf{Proof of Lemma 1}\\
From the construction of $C=(c^l)^\top \otimes I_d\in \mathbb{R}^{d\times md}$, where $c^l$ selects the $l$-th coordinate, we stress that the update (\ref{4e}) only involves the subvector $x_{t+1}^l$ as it can be seen by explicitly writing (\ref{4e}) as:
\begin{gather*}
    \lambda_{t+1} = \lambda_t +\tfrac{1}{\mu_2}(x_{t+1}^l-\theta_{t+1}).
\end{gather*}
We decompose the dual variable as $y^\top= \begin{bmatrix}\alpha^\top & \beta^\top \end{bmatrix}$, where $\alpha,\beta \in \mathbb{R}^{nd}$. Recalling $D^\top=\begin{bmatrix}I_{nd} & I_{nd} \end{bmatrix}$ and pre-multiplying the dual update (\ref{4d}) with $D^\top$, we obtain 
\begin{gather*}
    D^\top y_{t+1}=D^\top y_t+\tfrac{1}{\mu_1}D^\top (Ax_{t+1}-Dz_{t+1}).
\end{gather*}
Using (\ref{z}), we conclude that $D^\top y_{t+1}=0$ for $t \geq 0$, i.e., $\beta_{t+1}=-\alpha_{t+1}$. With zero initialization, we can therefore express $y_t$ as 
$
    y_t = \begin{bmatrix} 
    \alpha_t \\
    -\alpha_t
    \end{bmatrix}
$,
for all $t\geq 0$. This shows that it is redundant to maintain both $\alpha$ and $\beta$ since they sum to $0$ for all $t\geq 0$. Therefore, we can rewrite (\ref{4d}) as: 
\begin{gather}
    \alpha_{t+1}=\alpha_t+\tfrac{1}{\mu_1}(A_sx_{t+1}-z_{t+1}), \label{alpha1}\\
    -\alpha_{t+1} = -\alpha_{t}+\tfrac{1}{\mu_1}(A_dx_{t+1}-z_{t+1}). \label{alpha2}
\end{gather}
After taking the difference of the above two equations and dividing by 2, we obtain:
\begin{gather}
    \alpha_{t+1}=\alpha_t+\tfrac{1}{2\mu_1} E_sx_{t+1},
\end{gather}
where we have used the fact that $E_s=A_s-A_d$. Similarly, by summing (\ref{alpha1}) and (\ref{alpha2}), we obtain $z_{t+1}=\tfrac{1}{2}E_ux_{t+1}$. With zero initialization for $x_t$, it is redundant to keep $z_t$ since we can compute it as:
\begin{gather}
    z_t=\tfrac{1}{2}E_u x_{t} \,\,\forall\,\, t\geq 0. \label{z and x}
\end{gather}
Recalling the approximated Hessian in (\ref{approx Hessian}), we express the primal update (\ref{x}) as 
\begin{gather}
    x_{t+1}=x_t-H_t^{-1}\big\{\nabla F(x_t)
    +A^\top y_t+C^\top \lambda_t\nonumber \\ +\tfrac{1}{\mu_1}A^\top (Ax_t- Dz_t)+\tfrac{1}{\mu_2}C^\top (Cx_t-\theta_t) \big\}. \label{primal1}
\end{gather}
Since $E_s^\top = A_s^\top -A_d^\top$, $E_u^\top = A_s^\top +A_d^\top$, and $L_u= E_u^\top E_u$ as in (\ref{incidence}) and (\ref{laplacian}), along with the decomposition of $y_t$ and the equality in (\ref{z and x}), we have 
\begin{gather}
A^\top y_t=E_s^\top \alpha_t, \label{Esalpha}\\
\tfrac{1}{\mu_1}A^\top Dz_t = \tfrac{1}{2\mu_1}L_ux_t. \label{Luxt}
\end{gather}
Using (\ref{degree}), we can rewrite 
\begin{gather}
    \tfrac{1}{\mu_1}A^\top Ax_t-\tfrac{1}{2\mu_1}L_ux_t=\tfrac{1}{2\mu_1}(2\Delta-L_u)x_t=\tfrac{1}{2\mu_1}L_sx_t. \label{Lsxt}
\end{gather}
Substituting (\ref{Esalpha})-(\ref{Lsxt}) into (\ref{primal1}), we obtain the primal updates as in (\ref{14a}). \QEDB\\
\noindent \textbf{Proof of Lemma 2}\\
The update (\ref{14b}) is equivalent to 
\begin{gather*}
    \theta_{t+1} = \underset{\theta}{\mathrm{argmin}}\left\{g(\theta)+\tfrac{1}{2\mu_2}\norm{x_{t+1}^l+\mu_2\lambda_t-\theta}^2\right\}.
\end{gather*}
From the optimality condition, it holds that 
\begin{gather}
    0\in g(\theta_{t+1})-\tfrac{1}{\mu_2}(x_{t+1}+\mu_2\lambda_t-\theta_{t+1}).\label{plemma3a}
\end{gather}
Substituting the dual update $\lambda_{t+1}=\lambda_t+\tfrac{1}{\mu_2}(x_{t+1}^l-\theta_{t+1})$ into (\ref{plemma3a}), we obtain 
\begin{gather*}
    0\in \partial g(\theta_{t+1})-\lambda_{t+1}.
\end{gather*}
After re-arranging, we obtain the desired. \QEDB\\
\noindent\textbf{Proof of Lemma 3}\\
Note that the difference between BFGS-ADMM and it's first order variant lies in $H_t$ (\ref{approx Hessian}). For the first order variant, 
$
    H_t = \tfrac{1}{\mu_1}\Delta+\tfrac{1}{\mu_2}C^\top C+\epsilon I_{md},
$
i.e., $B_t$ is identically zero. By rearranging (\ref{14a}), we obtain the following equation:
\begin{gather}
    \nabla F(x_t)+E_s^\top \alpha_t+C^\top \lambda_t+\tfrac{1}{2\mu_1}L_sx_t+\tfrac{1}{\mu_2}C^\top (Cx_t-\theta_t)\nonumber\\
    +(B_t+\tfrac{1}{\mu_1}\Delta+\tfrac{1}{\mu_2}C^\top C+\epsilon I_{md})(x_{t+1}-x_t) = 0. \label{plemma2a}
\end{gather}
Using the dual update (\ref{14c}), we have
\begin{gather}
    E_s^\top \alpha_t+\tfrac{1}{2\mu_1}L_sx_t = E_s^\top\alpha_{t+1}-\tfrac{1}{2\mu_1}L_s(x_{t+1}-x_t). \label{plemma2b}
\end{gather}
Since $2\Delta=L_s+L_u$, it holds that 
\begin{gather}
    \tfrac{1}{\mu_1}\Delta(x_{t+1}-x_t)-\tfrac{1}{2\mu_1}L_s(x_{t+1}-x_t)=\tfrac{1}{2\mu_1}L_u(x_{t+1}-x_t).\label{plemma2c}
\end{gather}
Substituting (\ref{plemma2b}) and (\ref{plemma2c}) into (\ref{plemma2a}), we obtain
\begin{gather}
    \nabla F(x_t)+B_t(x_{t+1}-x_t)
    +\tfrac{1}{2\mu_1}(L_u+\epsilon I_{md})(x_{t+1}-x_t)\nonumber\\
    +C^\top \big\{\lambda_t+\tfrac{1}{\mu_2}(x_{t+1}^l-\theta_t)\big\}+E_s^\top\alpha_{t+1}=0. \label{plemma2e}
\end{gather}
Using the dual update (\ref{14d}), we have 
\begin{gather}
    \lambda_t+\tfrac{1}{\mu_2}(x_{t+1}^l-\theta_t) = \lambda_{t+1}+\tfrac{1}{\mu_2}(\theta_{t+1}-\theta_t). \label{plemma2d}
\end{gather}
Substituting (\ref{plemma2d}) into (\ref{plemma2e}) and using the definition of $e_t=\nabla F(x_t)+B_t(x_{t+1}-x_t)-\nabla F(x_{t+1})$, we obtain 
\begin{gather}
        \nabla F(x_{t+1})+\tfrac{1}{2\mu_1}(L_u+\epsilon I)(x_{t+1}-x_t)\nonumber\\
    +C^\top \big\{\lambda_{t+1}+\tfrac{1}{\mu_2}(\theta_{t+1}-\theta_t) \big\}+E_s^\top\alpha_{t+1}+e_t=0. \label{plemma2f}
\end{gather}
After subtracting (KKTa) from (\ref{plemma2f}), we obtain the desired where $B_t=0$ for the first order variant. \QEDB\\
\noindent\textbf{Proof of Lemma 4} \\
\noindent For BFGS-ADMM, since $B_{t+1}$ satisfies the secant equation: $B_{t+1}(x_{t+1}-x_t)=\nabla F(x_{t+1})-\nabla F(x_t)$, we have
\begin{gather*}
    e_{t}= (B_{t}-B_{t+1})(x_{t+1}-x_t),
\end{gather*}
and (\ref{bfgs error1}) follows by applying the Cauchy-Schwartz inequality. For the first order variant, the Lipschitz continuity of $\nabla F(x)$ implies that $\forall\,x,y \in\mathbb{R}^{md}$, 
$
    \norm{\nabla F(x)-\nabla F(y)} \leq M_f \norm{x-y}
$. The claim in (\ref{1st order error}) follows. \QEDB\\
\textbf{Proof of Lemma 5}\\
Consider the update (\ref{14c}) and (\ref{14d}). It can be written compactly as: 
$
    \begin{bmatrix}\alpha_{t+1} \\ \lambda_{t+1} \end{bmatrix} = \begin{bmatrix}\alpha_{t} \\ \lambda_{t} \end{bmatrix} +\begin{bmatrix}
    \tfrac{1}{2\mu_1}& 0 \\
    0& \tfrac{1}{\mu_2}
    \end{bmatrix}\begin{bmatrix} E_s \\ C \end{bmatrix}x_{t+1}-\begin{bmatrix} \textbf{0} \\\theta_{t+1} \end{bmatrix},
$
where $\textbf{0}$ is a zero vector of dimension $nd$. Therefore, it suffices to prove the vector $[\textbf{0}^\top \,\, \theta_{t+1}^\top]^\top$ lies in the column space of $X^\top$ . Note that since the graph is connected, $E_s r=0$ if $r\in \mathbb{R}^{md}$ is in consensus, i.e., $r^1=r^2=\dots=r^m \in \mathbb{R}^d$. By letting $r^i=\theta_{t+1}$ for all $i$, we have: 
$
    \begin{bmatrix} E_s \\ C \end{bmatrix}r= \begin{bmatrix} \textbf{0} \\\theta_{t+1} \end{bmatrix}
$,
which shows that $[\textbf{0}^\top \,\, \theta_{t+1}^\top]^\top$ lies in the column space of $X^\top$. To simplify the notation we denote $[\alpha_t^\top, \lambda_t^\top]^\top=w_t$. To prove the existence of optimal $w_\star$ in the column space of $X^\top$, consider the (KKTa) satisfied with any optimal $w$,
\begin{gather*}
    \nabla F(x_\star)+Xw=0.
\end{gather*}
The projection of $w$ into the column space of $X^\top$, denoted as $w_\star$, also satisfies (KKTa) since their differences lie in the kernel of $X$, i.e., $X(w-w_\star)=0$. The uniqueness of $w_\star$ can be established by contradiction. Let $w_1 = X^\top r_1$ and $w_2= X^\top r_2$ be two optimal stacked dual vector that lie in the column space of $X^\top$, $r_1\neq r_2$. Since $F(\cdot)$ is strongly convex, from (KKTa), it holds that 
$
    \nabla F(x_\star)+XX^\top r_1=0 ,\\
    \nabla F(x_\star)+XX^\top r_2=0.
$
After taking the difference, we have $XX^\top (r_1-r_2)=0$. But $XX^\top = E_s^\top E_s+C^\top C=L_s+C^\top C>0$. Therefore $r_1-r_2=0$, contradiction. Inequality (\ref{lemma5}) follows from the fact that $w_{t+1}-w_\star$ is orthogonal to the kernel of $X$.\QEDB

\end{document}